\documentclass[a4paper,oneside]{amsart}
\usepackage[german,english]{babel}

\RequirePackage{amsmath}
\RequirePackage{bm}
\RequirePackage{amssymb}
\RequirePackage{upref}
\RequirePackage{amsthm}
\RequirePackage{enumerate}
\RequirePackage{pb-diagram}
\RequirePackage{amsfonts}
\RequirePackage[mathscr]{eucal}
\RequirePackage{verbatim}
\RequirePackage{xr}
\RequirePackage[pdftex]{graphicx}
\RequirePackage[pdftex]{floatflt}

\newcommand{\al}{\alpha}
\newcommand{\bet}{\beta}
\newcommand{\ga}{\gamma}
\newcommand{\de}{\delta }
\newcommand{\e}{\epsilon}

\newcommand{\f}{\varphi}
\newcommand{\h}{\eta}

\newcommand{\ka}{\kappa}
\newcommand{\lam}{\lambda}
\newcommand{\m}{\mu}
\newcommand{\n}{\nu}

\newcommand{\s}{\sigma}
\newcommand{\x}{\xi}

\newcommand{\C}{\varGamma}

\newcommand{\F}{\varPhi}

\newcommand{\Om}{\varOmega}

\newcommand{\Si}{\varSigma}

\newcommand{\di}[1]{#1\nobreakdash-\hspace{0pt}dimensional}%\di n
\newcommand{\nbdd}{\nobreakdash--}
\newcommand{\nbd}{\nobreakdash-\hspace{0pt}}

\newcommand{\fm}[1]{F_{|_{M_#1}}}
\newcommand{\fmo}[1]{F_{|_{#1}}}%\fmo M

\newcommand{\fv}[2]{#1\hspace{0pt}_{|_{#2}}}

\newcommand{\so}{{\mc S_0}}

\newcommand{\const}{\tup{const}}

\newcommand{\msp[1]}[1]{\mspace{#1mu}}
\newcommand{\low}[1]{{\hbox{}_{#1}}}

\newcommand{\R}[1][n+1]{{\protect\mathbb R}^{#1}}

\newcommand{\N}{{\protect\mathbb N}}

\newcommand{\eR}{\stackrel{\lower1ex \hbox{\rule{6.5pt}{0.5pt}}}{\msp[3]\R[]}}
\newcommand{\eN}{\stackrel{\lower1ex \hbox{\rule{6.5pt}{0.5pt}}}{\msp[1]\N}}
\newcommand{\eO}{\stackrel{\lower1ex
\hbox{\rule{6pt}{0.5pt}}}{\msc O}}

\DeclareMathOperator{\graph}{graph}

\newcommand\pa{\partial}
\newcommand\pde[2]{\frac {\partial#1}{\partial#2}}
\newcommand\pd[3]{\frac {\partial#1}{\partial#2^#3}}   %e.g. \pd fxi
\newcommand\pdc[3]{\frac {\partial#1}{\partial#2_#3}}   %contravariant
\newcommand\pdm[4]{\frac {\partial#1}{\partial#2_#3^#4}}   %mixed
\newcommand\pddc[4]{\frac {{\partial\hskip0.15em}^2#1}{\partial {#2_
#3}\,\partial{#2_#4}}} 
\newcommand{\un}{\infty}
\newcommand{\A}{\forall}

\newcommand{\set}[2]{\{\,#1\colon #2\,\}}
\newcommand{\uu}{\cup}
\newcommand{\ii}{\cap}
\newcommand{\uuu}{\bigcup}

\newcommand{\uud}{ \stackrel{\lower 1ex \hbox {.}}{\uu}}
\newcommand{\uuud}[1]{ \stackrel{\lower 1ex \hbox {.}}{\uuu_{#1}}}
\newcommand\su{\subset}

\newcommand{\sminus}[1][28]{\raise 0.#1ex\hbox{$\scriptstyle\setminus$}}

\newcommand{\abs}[1]{\lvert#1\rvert}

\newcommand{\norm}[1]{\lVert#1\rVert}

\newcommand{\nnorm}[1]{| \mspace{-2mu} |\mspace{-2mu}|#1| \mspace{-2mu}
|\mspace{-2mu}|}
\newcommand{\spd}[2]{\protect\langle #1,#2\protect\rangle}

\newcommand\ch[3]{\varGamma_{#1#2}^#3}
\newcommand\cha[3]{{\bar\varGamma}_{#1#2}^#3}

\newcommand{\riem}[4]{R_{#1#2#3#4}}
\newcommand{\riema}[4]{{\bar R}_{#1#2#3#4}}

\newcommand{\tbf}{\textbf}
\newcommand{\tit}{\textit}

\newcommand{\tup}{\textup}% text upright

\newcommand{\mc}{\protect\mathcal}
\newcommand{\msc}{\protect\mathscr}

\providecommand{\bysame}{\makebox[3em]{\hrulefill}\thinspace}

\newcommand{\ci}{\cite}
\newcommand{\bib}{\bibitem}

\newcommand{\bt}{\begin{thm}}
\newcommand{\bl}{\begin{lem}}
\newcommand{\bc}{\begin{cor}}
\newcommand{\bd}{\begin{definition}}
\newcommand{\bpp}{\begin{prop}}
\newcommand{\br}{\begin{rem}}
\newcommand{\bn}{\begin{note}}
\newcommand{\be}{\begin{ex}}
\newcommand{\bes}{\begin{exs}}
\newcommand{\bb}{\begin{example}}
\newcommand{\bbs}{\begin{examples}}
\newcommand{\ba}{\begin{axiom}}

\newcommand{\et}{\end{thm}}
\newcommand{\el}{\end{lem}}
\newcommand{\ec}{\end{cor}}
\newcommand{\ed}{\end{definition}}
\newcommand{\epp}{\end{prop}}
\newcommand{\er}{\end{rem}}
\newcommand{\en}{\end{note}}
\newcommand{\ee}{\end{ex}}
\newcommand{\ees}{\end{exs}}
\newcommand{\eb}{\end{example}}
\newcommand{\ebs}{\end{examples}}
\newcommand{\ea}{\end{axiom}}

\newcommand{\bp}{\begin{proof}}
\newcommand{\ep}{\end{proof}}
\newcommand{\eps}{\renewcommand{\qed}{}\end{proof}}

\newcommand{\bal}{\begin{align}}

\newcommand{\bi}[1][1.]{\begin{enumerate}[\upshape #1]}
\newcommand{\bia}[1][(1)]{\begin{enumerate}[\upshape #1]}
\newcommand{\bin}[1][1]{\begin{enumerate}[\upshape\bfseries #1]}
\newcommand{\bir}[1][(i)]{\begin{enumerate}[\upshape #1]}
\newcommand{\bic}[1][(i)]{\begin{enumerate}[\upshape\hspace{2\cma}#1]}
\newcommand{\bis}[2][1.]{\begin{enumerate}[\upshape\hspace{#2\parindent}#1]}
\newcommand{\ei}{\end{enumerate}}

\newcommand\ndots{\raise 0.47ex \hbox {,}\hskip0.06em\cdots %
     \raise 0.47ex \hbox {,}\hskip0.06em} 

\newcommand{\q}{\quad}
\newcommand{\qq}{\qquad}

\newcommand{\hp}{\hphantom}

\newcommand\nd{\noindent}

\newskip\Csmallskipamount                                                
\Csmallskipamount=\smallskipamount
\newskip\Cmedskipamount
\Cmedskipamount=\medskipamount
\newskip\Cbigskipamount
\Cbigskipamount=\bigskipamount

\newcommand\cvs{\vspace\Csmallskipamount}   
\newcommand\cvm{\vspace\Cmedskipamount}
\newcommand\cvb{\vspace\Cbigskipamount}

\newskip\csa
\csa=\smallskipamount

\newskip\cma
\cma=\medskipamount

\newskip\cba
\cba=\bigskipamount

\newdimen\spt
\newcommand\citem{\cvs\advance\itemno by
1{(\romannumeral\the\itemno})\hskip3pt}
\newcommand{\bitem}{\cvm\nd\advance\itemno by
1{\bf\the\itemno}\hspace{\cma}}

\newcount\itemno
\newcommand{\las}[1]{\label{S:#1}}

\newcommand{\lae}[1]{\label{E:#1}}
\newcommand{\lat}[1]{\label{T:#1}}
\newcommand{\lal}[1]{\label{L:#1}}
\newcommand{\lad}[1]{\label{D:#1}}

\newcommand{\lar}[1]{\label{R:#1}}

\newcommand{\rs}[1]{Section~\ref{S:#1}}

\newcommand{\rl}[1]{Lemma~\ref{L:#1}}

\newcommand{\rr}[1]{Remark~\ref{R:#1}}

\newcommand{\re}[1]{\eqref{E:#1}}

\RequirePackage{upref}
\RequirePackage{amsthm}
\RequirePackage{enumerate}%\begin{enumerate}[(i)]

\theoremstyle{plain}
\newtheorem{thm}{Theorem}[section]
\newtheorem{lem}[thm]{Lemma}
\newtheorem{prop}[thm]{Proposition}
\newtheorem{cor}[thm]{Corollary}

\theoremstyle{definition}
\newtheorem{rem}[thm]{Remark}
\newtheorem{definition}[thm]{Definition}
\newtheorem{example}[thm]{Example}
\newtheorem{ex}[thm]{Exercise}

\swapnumbers
\theoremstyle{remark}

\numberwithin{equation}{section}
\usepackage{xr-hyper}
\usepackage[hyperindex=true, pdfauthor= Claus\  Gerhardt, pdftitle= Closed\ Weingarten\ hypersurfaces, bookmarks=true, extension= pdf, colorlinks=false, plainpages=false,hyperfootnotes=false, debug=false, pagebackref]{hyperref}

\newcommand{\anl}{\htmladdnormallink}

\begin{document}
%\larger[1]
\title[Closed Weingarten hypersurfaces]{Closed Weingarten hypersurfaces in
semi-riemannian manifolds}

% author one information
\author{Claus Gerhardt}
\address{Ruprecht-Karls-Universit\"at, Institut f\"ur Angewandte Mathematik,
Im Neuenheimer Feld 294, 69120 Heidelberg, Germany}
%\curraddr{}
\email{gerhardt@math.uni-heidelberg.de}
\urladdr{http://www.math.uni-heidelberg.de/studinfo/gerhardt/}
%\thanks{}

% author two information
%\author{}
%\address{}
%\curraddr{}
%\email{}
%\thanks{}
%
\date{November 6, 2001}
%
% at present the "communicated by" line appears only in ERA and PROC
%\commby{}

%\dedicatory{}

\begin{abstract} The existence of closed hypersurfaces of prescribed curvature
in semi-riemannian  manifolds is proved provided there are barriers.
\end{abstract}
\maketitle
\thispagestyle{empty}

\tableofcontents
\setcounter{section}{-1}
\section{Introduction} 

\cvb
We want to prove the existence of closed hypersurfaces of prescribed
curvature in Riemannian or Lorentzian manifolds $N$ of dimension $n~+~1$, $n\ge
2$. Since we wish to treat both cases simultaneously, let us stipulate that terms
which make only sense in Lorentzian manifolds should be ignored if the ambient
space is Riemannian. With this in mind, let $\Om$ be a connected, precompact, 
open subset of
$N$, $f=f(x,\n)$ a positive function defined for $x\in\bar\Om$ and time-like
vectors $\n\in T_x(N)$, and $F\in C^{2,\al}(\C_+)\ii C^0(\bar\C_+)$ a symmetric
curvature function defined on the positive cone $\C_+\su \R[n]$. Then, we look
for closed space-like hypersurfaces $M\su\Om$ such that

\begin{equation}\lae{0.1}
\fv F M=f(x,\n),
\end{equation}

\cvm
\nd where $f$ is evaluated at $x\in M$ and at the past directed normal $\n\in
T_x(N)$.

\cvm
Various existence results have been proved for a wide range of curvature
functions $F$ if $f$ only depends on $x$. In Euclidian space any monotone
curvature function $F$ can be considered with the property that $\log F$ is
concave---at least in principle, cf. \ci{cns1, cg2}. The only possible obstruction
could occur when one tries to prove $C^1$- estimates. That particular difficulty
arises in any Riemannian space, but in addition---in order to obtain $C^2$-
estimates---one has to assume that $F$ satisfies a certain concavity estimate,
i.e. a stronger property is needed than mere concavity of $F$ or $\log F$. In
\ci{cg2,  cg4} we proved existence results for curvature functions of class
$(K)$ that satisfy such an estimate.

\cvm
If the ambient space is Lorentzian, then, $C^1$- estimates can be obtained for
curvature functions for which a corresponding estimate in a Riemannian setting
is known to be impossible, or only achievable with additional structural
conditions on $f$. But still, the curvature functions have to satisfy the same
stronger concavity property as in the Riemannian case in order to derive $C^2$-
estimates, and, in addition, one further estimate is needed, namely, there should
exist
$\e_0>0$ such that

\begin{equation}
\e_0 F H\le F^{ij}h_{ik}h^k_j
\end{equation}

\cvm
\nd for all admissible tensors $(h_{ij})$, where as usual $H$ stands for the mean
curvature, cf. \ci{cg8}.

\cvm
Recently, a concavity estimate has been proved for the scalar curvature
operator $H_2$, cf. \ci{pb2}. We  improved that estimate in \ci{cg9}, so that
functions $f$ depending on the normal can be considered, and had been able to
prove the existence of closed space-like hypersurfaces satisfying \re{0.1},
where $N$ is Lorentzian, $F=H_2$, and $f=f(x,\n)$ is general enough, so that
solutions can be considered as having prescribed scalar curvature.

\cvm
It is only natural to ask if similar generalizations with regard to $f$ are also
possible for other curvature functions or in Riemannian spaces instead of
Lorentzian. Two difficulties arise from the presence of the normal vector in the
right-hand side $f$ that have to be dealt with separately. Let us first address
the simple one, the $C^2$- estimates. These estimates can be derived for all
curvature functions $F$ that obey a strong concavity condition mentioned above,
and, in addition, an estimate of the form

\begin{equation}\lae{0.3}
n\e_0 F\le F_i H\qq\A\,1\le i\le n,
\end{equation}

\cvm
\nd where $\e_0=\e_0(F)$ is a positive constant, $F_i=\pde F{\ka_i}$, and the
inequality should be valid in the convex cone associated with $F$. Furthermore,
it is assumed that lower order estimates in the $C^0$- and $C^1$- norms are
already established.

\cvm
Obtaining the $C^1$- estimates is the second difficulty that seems to be
unsurmountable in general Riemannian spaces, e.g. for $H_2$, and even in
Euclidean space, it is only possible if special structural conditions on $f$ are
imposed. However, the only curvature functions $F$ for which strong concavity
estimates are known so far, are $H_2$ and those of class $(K)$. The latter are
defined in $\C_+$, i.e. the admissible hypersurfaces have to be strictly convex,
and, thus, the $C^1$- estimates are easily obtained.

\cvm
We have excluded the mean curvature or functions of mean curvature type, $F\in
(H)$, from the solvability discussion, since they pose different problems. In
their case the $C^1$- estimates are the only challenge, and they can be derived
in Lorentzian space under the assumptions

\begin{align}
\nnorm{f_\bet(x,\n)}&\le c (1+\nnorm{\n}),\\
\intertext{and}
\nnorm{f_{\n^\bet}(x,\n)}&\le c,
\end{align}

\cvm
\nd cf. \ci[Proposition 4.8]{cg9}.

\cvb
To give a precise statement of the existence results we need a few definitions
and assumptions. It seems advisable to treat the Lorentzian and Riemannian
cases separately.

\subsection{The Lorentzian case} We assume that $N$ is a smooth, connected,
\tit{globally hyperbolic} manifold with a compact \tit{Cauchy hypersurface}
$\so$,  and suppose that $\Om$ is bounded by  two {\it achronal}, connected,
space-like hypersurfaces $M_1$ and $M_2$ of class $C^{4,\al}$, where $M_1$ is
supposed to lie in the past of $M_2$. 

\cvm
Let $F$  of class $(K)$ satisfy (0.3), and $0<\e_1\le f$ be of class $C^{2,\al}$. Then,
we assume that the boundary components act as barriers for $(F,f)$.

\cvm
\bd\lad{0.1}
$M_2$ is an \tit{upper barrier} for $(F,f)$, if $M_2$ is strictly convex and
satisfies
\begin{equation}
\fm 2\ge f(x,\n),
\end{equation}
where $f$ is evaluated at $x\in M_2$, and at the past directed normal $\n(x)$ of
$M_2$.

$M_1$ is a
\tit{lower barrier} for $(F,f)$, if at the points
$\Si\su M_1$, where
$M_1$ is strictly convex, there holds

\begin{equation}
\fv F\Si\le f(x,\n).
\end{equation}

\cvm\nd
$\Si$ may be empty.
\ed

\cvm
Then, we can prove

\bt\lat{0.2}
Let $M_1$ be a lower and $M_2$ an upper barrier for $(F,f)$. Then, the problem

\begin{equation}
\fmo M= f(x,\n)
\end{equation}

\cvm\nd
has a strictly convex solution $M\su \bar\Om$ of class $C^{4,\al}$ that can be
written as a graph over $\mc S_0$ provided there exists a strictly convex
function $\chi\in C^2(\bar\Om)$.
\et

\cvm
\subsection{The Riemannian case} Let $N$ be a smooth connected Riemannian
manifold with $K_N\le 0$, and assume that the boundary components of $\Om$ are
both strictly convex hypersurfaces homeomorphic to $S^n$ and of class
$C^{4,\al}$, such that the mean curvature vector of $M_1$ points outside of $\Om$
and the mean curvature vector of $M_2$ points inside of $\Om$.

\cvm
Then, we can prove

\cvm
\bt\lat{0.3}
Let $F$  of class $(K)$ satisfy (0.3), and let $0<\e_1\le f$ of class $C^{2,\al}$ be
given. Assume that $M_2$ is an upper barrier for $(F,f)$ and $M_1$ a lower barrier.
Then,  the problem

\begin{equation}
\fmo M= f(x,\n)
\end{equation}

\cvm\nd
has a strictly convex solution $M\su \bar\Om$ of class $C^{4,\al}$. 
\et

\cvb
The paper is organized as follows: In \rs 1 we take a closer look at   curvature
functions and show that the curvature functions in question coincide with a
subclass of $(K)$, where the functions can be written as a product such that
one factor is a power of the Gaussian curvature.  

\cvm

In \rs 2 we introduce the notations and common definitions we rely on,  and state
the equations of Gau{\ss}, Codazzi, and Weingarten for space-like hypersurfaces
in a semi-riemannian manifold.

\cvm
In \rs 3 we look at the curvature flow associated with our problem, and the
corresponding evolution equations for the basic geometrical quantities of the
flow hypersurfaces.

\cvm
In \rs 4 we prove lower order estimates for the evolution problem,  while a
priori estimates in the
$C^2$\nbd{norm} are derived in
\rs 5 for the Lorentzian case, and in \rs 6 for the Riemannian case.

\cvm
The final existence result is contained in \rs 7.

\cvb
\section{Curvature functions}\las{1}

Let $\C_+\su\R[n]$ be the open positive cone and $F\in C^{2,\al}(\C_+)\ii
C^0(\bar\C_+)$ a symmetric function satisfying the condition

\begin{equation}\lae{1.1}
F_i=\pd F\ka i>0\; ;
\end{equation}

\cvm
\nd
then, $F$ can also be viewed as a function defined on the space of symmetric,
positive definite matrices $\mathscr S_+$, for, let $(h_{ij})\in \msc S_+$ with
eigenvalues $\ka_i,\,1\le i\le n$, then define $F$ on $\msc S_+$ by

\begin{equation}
F(h_{ij})=F(\ka_i).
\end{equation}

If we define 
\begin{align}
F^{ij}&=\pde F{h_{ij}}\\
\intertext{and}
F^{ij,kl}&=\pddc Fh{{ij}}{{kl}}
\end{align}
then, 
\begin{equation}
F^{ij}\x_i\x_j=\pdc F\ka i \abs{\x^i}^2\q\A\, \x\in\R[n],
\end{equation}
\begin{equation}
F^{ij} \,\text{is diagonal if $h_{ij}$ is diagonal,}
\end{equation}
and
\begin{equation}\lae{1.7}
F^{ij,kl}\h_{ij}\h_{kl}=\pddc F{\ka}ij\h_{ii}\h_{jj}+\sum_{i\ne
j}\frac{F_i-F_j}{\ka_i-\ka_j}(\h_{ij})^2,
\end{equation}
for any $(\h_{ij})\in \msc S$, where $\msc S$ is the space of all symmetric
matrices. The second term on the right-hand side of \re{1.7} is non-positive if
$F$ is concave, and non-negative if $F$ is convex, and has to be interpreted as a
limit if $\ka_i=\ka_j$.

\cvm
In \ci{cg8} we defined---or better redefined---the curvature functions of class
$(K)$ as

\cvm
\bd
A symmetric curvature function $F\in C^{2,\al}(\C_+)\ii C^0(\bar \C_+)$
positively homogeneous of degree $d_0>0$ is said to be of class $(K)$ if

\begin{equation}\lae{1.8}
F_i=\pd F{\ka}i>0\q \text{in } \C_+,
\end{equation}
\begin{equation}\lae{1.9}
\fmo{\pa \C_+}=0,
\end{equation}
and
\begin{equation}\lae{1.10}
F^{ij,kl}\h_{ij}\h_{kl}\le F^{-1}(F^{ij}\h_{ij})^2-F^{ik}\tilde
h^{jl}\h_{ij}\h_{kl}\qq\A\,\h\in\msc S,
\end{equation}

\cvm
\nd
or, equivalently, if we set $\hat F=\log F$,

\begin{equation}\lae{1.11}
\hat F^{ij,kl}\h_{ij}\h_{kl}\le -\hat F^{ik}\tilde
h^{jl}\h_{ij}\h_{kl}\qq\A\,\h\in\msc S,
\end{equation}

\cvm
\nd
where $F$ is evaluated at $(h_{ij})$.
\ed

\cvb
The preceding considerations are
also applicable if the
$\ka_i$ are the principal curvatures of a hypersurface $M$ with metric $(g_{ij})$.
$F$ can then be looked at as being defined on the space of all symmetric tensors
$(h_{ij})$ with eigenvalues $\ka_i$ with respect to the metric.

\begin{equation}
F^{ij}=\pdc Fh{{ij}}
\end{equation}

\cvm
\nd
is then a contravariant tensor of second order. Sometimes it will be convenient
to circumvent the dependence on the metric by considering $F$ to depend on the
mixed tensor

\begin{equation}
h_j^i=g^{ik}h_{kj}.
\end{equation}

\cvm\nd
Then,

\begin{equation}
F_i^j=\pdm Fhji
\end{equation}

\cvm
\nd
is also a mixed tensor with contravariant index $j$ and covariant index $i$.

\cvb
\br
Let $F\in (K)$, then $\log F$ is concave, and, if $F$ is homogeneous of degree $1$,
then, $F$ is already concave.
\er

\bp
The concavity of $\log F$ follows immediately from \re{1.11}, while, in case $F$ is
homogeneous of degree $1$, the concavity of $F$ can be derived from the
inequality \re{1.10} by applying Schwartz inequality: Choose coordinates such
that in a fixed point

\begin{equation}\lae{1.15}
g_{ij}=\de_{ij}\q\tup{and}\q h_{ij}=\ka_i\de_{ij}.
\end{equation}

\cvm
Let $\h=(\h_{ij})$ be an arbitrary symmetric tensor, then

\begin{equation}
\begin{aligned}
F^{-1}(F^{ij}\h_{ij})^2&=F^{-1}\big(\sum_iF^i_i\h_{ii}\big)^2\\[\cma]
&=
F^{-1}\big(\sum_i(F^i_i)^\frac12\ka_i^\frac12(F^i_i)^\frac12\ka_i^{-\frac12}
\h_{ii}\big)^2\\[\cma]
&\le F^{-1}\big(\sum_iF^i_i\ka_i\big)\msp \big(\sum_iF^i_i\ka_i^{-1}\h_{ii}^2\big)
\\[\cma]
&\le F^{ik} \tilde h^{jl}\h_{ij}\h_{kl},
\end{aligned}
\end{equation}

\cvm
\nd hence, the right-hand side of \re{1.10} is non-positive.
\ep

\cvm
The subclass $(K^*)$ has been defined in \ci[Definition 1.6]{cg8} as

\cvm
\bd\lad{1.3}
A function $F\in (K)$ is said to be of class $(K^*)$ if there exists
$0<\e_0=\e_0(F)$ such that

\begin{equation}\lae{1.17}
\e_0 F H\le F^{ij} h_{ik}h^k_j\,,
\end{equation}

\cvm
\nd
for any $(h_{ij})\in \msc S_+$, where $F$ is evaluated at $(h_{ij})$. $H$ represents
the mean curvature, i.e. the trace of $(h_{ij})$.
\ed

\cvm
The condition \re{1.17} is crucial for solving curvature problems in Lorentzian
manifolds; it is slightly weaker than the condition \re{0.3} which is also satisfied
by $F=H_2$, see e.g. \ci[(1.17)]{cg9}.

\cvm
\bl
Let $F$ be a symmetric curvature function defined in an  open convex cone $\C$
satisfying the  relations \re{1.1} and \re{0.3}
in
$\C$, then it also satisfies the inequality \re{1.17} with the same constant $\e_0$.
\el

\bp
We first observe that in view of the relation \re{0.3} $H$ is positive in
$\C$. Next, choose coordinates as in
\re{1.15}, then,

\begin{equation}
\begin{aligned}
F^{ij}h_{ik}h^k_j&=H^{-1} \sum_iH F_i\ka_i^2\\[\cma]
&\ge n\e_0 H^{-1}F\abs A^2\ge \e_0 F H,
\end{aligned}
\end{equation}

\cvm
\nd where we used the usual abbreviation $\abs A^2$ for $\sum_i\ka_i^2$.
\ep

\cvm
\br
Special functions of class $(K^*)$ are those that can be written as a product

\begin{equation}
F=G K^a,\qq a>0,
\end{equation}

\cvm
\nd where $G\in (K)$ and $K$ is the Gaussian curvature. They are exactly those
that satisfy the estimate

\begin{equation}\lae{1.20}
F_i\ka_i\ge \e_0 F\qq\A\,1\le i\le n
\end{equation}

\cvm
\nd with some positive constant $\e_0=\e_0(F)$, cf. \ci[Proposition 1.9]{cg8}.

\cvm
Using the simple estimate $\ka_i\le H$, which is valid in $\C_+$, we conclude that
these special functions also satisfy the condition \re{0.3}. 
\er

\cvm
The reverse is also
true.

\bl
Let $F\in (K)$ be such that the relation \re{0.3} is valid, then, $F$ also satisfies
\re{1.20}.
\el

\bp
Let $\ka_n$ be the largest component of the $n$- tupel $(\ka_i)\in \C_+$. Then, we
conclude in view of
\re{0.3} 

\begin{equation}
F_n\ka_n\ge \tfrac1nF_nH\ge \e_0 F.
\end{equation}

\cvm
On the other hand, any curvature function of class $(K)$ satisfies

\begin{equation}
F_i\ka_i\ge F_n\ka_n\qq\A\,1\le i\le n,
\end{equation}

\cvm\nd
cf.  \ci[Lemma 1.3]{cg4}.
\ep

\cvb
\section{Notations and preliminary results}\las 2

The main objective of this section is to state the equations of Gau{\ss}, Codazzi,
and Weingarten for hypersurfaces. We shall formulate the governing equations of
a hypersurface $M$ in a semi-riemannian \di{(n+1)} space $N$, which is either
Riemannian or Lorentzian. Geometric quantities in $N$ will be denoted by
$(\bar g_{\al\bet}),(\riema \al\bet\ga\de)$, etc., and those in $M$ by $(g_{ij}), (\riem
ijkl)$, etc. Greek indices range from $0$ to $n$ and Latin from $1$ to $n$; the
summation convention is always used. Generic coordinate systems in $N$ resp.
$M$ will be denoted by $(x^\al)$ resp. $(\x^i)$. Covariant differentiation will
simply be indicated by indices, only in case of possible ambiguity they will be
preceded by a semicolon, i.e. for a function $u$ in $N$, $(u_\al)$ will be the
gradient and
$(u_{\al\bet})$ the Hessian, but e.g., the covariant derivative of the curvature
tensor will be abbreviated by $\riema \al\bet\ga{\de;\e}$. We also point out that

\begin{equation}
\riema \al\bet\ga{\de;i}=\riema \al\bet\ga{\de;\e}x_i^\e
\end{equation}

\cvm
\nd
with obvious generalizations to other quantities.

\cvm
Let $M$ be a \tit{space-like} hypersurface, i.e. the induced metric is Riemannian,
with a differentiable normal $\n$. We define the signature of $\n$, $\s=\s(\n)$, by

\begin{equation}
\s=\bar g_{\al\bet}\n^\al\n^\bet=\spd \n\n.
\end{equation}

\cvm
In case $N$ is Lorentzian, $\s=-1$, and $\n$ is time-like.

\cvm
In local coordinates, $(x^\al)$ and $(\x^i)$, the geometric quantities of the
space-like hypersurface $M$ are connected through the following equations

\begin{equation}\lae{2.3}
x_{ij}^\al=-\s h_{ij}\n^\al
\end{equation}

\cvm
\nd
the so-called \tit{Gau{\ss} formula}. Here, and also in the sequel, a covariant
derivative is always a \tit{full} tensor, i.e.

\begin{equation}
x_{ij}^\al=x_{,ij}^\al-\ch ijk x_k^\al+\cha \bet\ga\al x_i^\bet x_j^\ga.
\end{equation}

\cvm
\nd
The comma indicates ordinary partial derivatives.

In this implicit definition the \tit{second fundamental form} $(h_{ij})$ is taken
with respect to $-\s\n$.

The second equation is the \tit{Weingarten equation}

\begin{equation}
\n_i^\al=h_i^k x_k^\al,
\end{equation}

\cvm\nd
where we remember that $\n_i^\al$ is a full tensor.

\cvm
Finally, we have the \tit{Codazzi equation}

\begin{equation}
h_{ij;k}-h_{ik;j}=\riema\al\bet\ga\de\n^\al x_i^\bet x_j^\ga x_k^\de
\end{equation}
and the \tit{Gau{\ss} equation}
\begin{equation}
\riem ijkl=\s \{h_{ik}h_{jl}-h_{il}h_{jk}\} + \riema \al\bet\ga\de x_i^\al x_j^\bet x_k^\ga
x_l^\de.
\end{equation}

\cvb
For the rest of this section we treat the Riemannian and Lorentzian cases
separately.

\subsection{The Lorentzian case}

Now, let us assume that $N$ is a globally hyperbolic Lorentzian manifold with a
\tit{compact} Cauchy surface $\so$.
Then, $N$ is topologically a product, $N=\R[]\times\mc{S}_0$, where $\mc S_0$ is
a compact, \di n Riemannian manifold, and there exists a Gaussian coordinate
system $(x^\al)_{0\le \al\le n}$ such that $x^0$ represents the time, the
$(x^i)_{1\le i\le n}$ are local coordinates for $\mc S_0$, where we may assume
that $\mc S_0$ is equal to the level hypersurface $\{x^0=0\}$---we don't
distinguish between $\mc S_0$ and $\{0\}\times \mc S_0$---, and such that the
Lorentzian metric takes the form

\begin{equation}\lae{2.8}
d\bar s_N^2=e^{2\psi}\{-{dx^0}^2+\s_{ij}(x^0,x)dx^idx^j\},
\end{equation}

\cvm\nd
where $\s_{ij}$ is a Riemannian metric, $\psi$ a function on $N$, and $x$ an
abbreviation for the space-like components $(x^i)$, see \ci{GR},
\ci[p.~212]{HE}, \ci[p.~252]{GRH}, and \ci[Section~6]{cg1}.
 We also assume that
the coordinate system is \tit{future oriented}, i.e. the time coordinate $x^0$
increases on future directed curves. Hence, the \tit{contravariant} time-like
vector $(\x^\al)=(1,0,\dotsc,0)$ is future directed as is its \tit{covariant}
version
$(\x_\al)=e^{2\psi}(-1,0,\dotsc,0)$.

\cvm
Furthermore, any achronal hypersurface can be written as a graph over $\so$, cf.
\ci[Proposition 2.5]{cg8}.

\cvm
 Let $M=\graph \fv u\so$ be a space-like hypersurface

\begin{equation}
M=\set{(x^0,x)}{x^0=u(x),\,x\in\mc S_0},
\end{equation}

\cvm
\nd
then the induced metric has the form

\begin{equation}
g_{ij}=e^{2\psi}\{-u_iu_j+\s_{ij}\}
\end{equation}

\cvm\nd
where $\s_{ij}$ is evaluated at $(u,x)$, and its inverse $(g^{ij})=(g_{ij})^{-1}$ can
be expressed as

\begin{equation}\lae{2.10}
g^{ij}=e^{-2\psi}\{\s^{ij}+\frac{u^i}{v}\frac{u^j}{v}\},
\end{equation}

\cvm\nd
where $(\s^{ij})=(\s_{ij})^{-1}$ and

\begin{equation}\lae{2.11}
\begin{aligned}
u^i&=\s^{ij}u_j\\
v^2&=1-\s^{ij}u_iu_j\equiv 1-\abs{Du}^2.
\end{aligned}
\end{equation}

\cvm
Hence, $\graph u$ is space-like if and only if $\abs{Du}<1$.

\cvm
The covariant form of a normal vector of a graph looks like

\begin{equation}
(\n_\al)=\pm v^{-1}e^{\psi}(1, -u_i).
\end{equation}

\cvm\nd
and the contravariant version is

\begin{equation}
(\n^\al)=\mp v^{-1}e^{-\psi}(1, u^i).
\end{equation}

\cvm
Thus, we have

\cvm
\br Let $M$ be space-like graph in a future oriented coordinate system. Then, the
contravariant future directed normal vector has the form

\begin{equation}
(\n^\al)=v^{-1}e^{-\psi}(1, u^i)
\end{equation}
and the past directed
\begin{equation}\lae{2.15}
(\n^\al)=-v^{-1}e^{-\psi}(1, u^i).
\end{equation}
\er

\cvm
In the Gau{\ss} formula \re{2.3} we are free to choose the future or past directed
normal, but we stipulate that we always use the past directed normal for reasons
that we have explained in \ci[Section 2]{cg8}.

\cvm
Look at the component $\al=0$ in \re{2.3} and obtain in view of \re{2.15}

\begin{equation}\lae{2.16}
e^{-\psi}v^{-1}h_{ij}=-u_{ij}-\cha 000\mspace{1mu}u_iu_j-\cha 0i0
\mspace{1mu}u_j-\cha 0j0\mspace{1mu}u_i-\cha ij0.
\end{equation}
Here, the covariant derivatives a taken with respect to the induced metric of
$M$, and
\begin{equation}
-\cha ij0=e^{-\psi}\bar h_{ij},
\end{equation}
where $(\bar h_{ij})$ is the second fundamental form of the hypersurfaces
$\{x^0=\const\}$.

\cvm
Sometimes, we need a Riemannian reference metric, e.g. if we want to estimate
tensors. Since the Lorentzian metric can be expressed as

\begin{equation}
\bar g_{\al\bet}dx^\al dx^\bet=e^{2\psi}\{-{dx^0}^2+\s_{ij}dx^i dx^j\},
\end{equation}

\cvm\nd
we define a Riemannian reference metric $(\tilde g_{\al\bet})$ by

\begin{equation}
\tilde g_{\al\bet}dx^\al dx^\bet=e^{2\psi}\{{dx^0}^2+\s_{ij}dx^i dx^j\}
\end{equation}

\cvm\nd
and we abbreviate the corresponding norm of a vectorfield $\h$ by

\begin{equation}
\nnorm \h=(\tilde g_{\al\bet}\h^\al\h^\bet)^{1/2},
\end{equation}

\cvm\nd
with similar notations for higher order tensors.

\cvb
\subsection{The Riemannian case} In view of our assumptions on $N$ and $\Om$, we
may assume that $N$ is simply connected and that $\Om$ is the difference of two
convex bodies, cf. \ci[Theorem 4.7]{cg2}, and, therefore, $\Om$ can be covered
by a geodesic polar coordinate system $(x^\al)_{0\le\al\le n}$, where $x^0$
re\-pre\-sents the radial distance to the center, and $(x^i)$ are local
coordinates for the geodesic sphere $S^n=\{x^0=1\}$.

\cvm
The barriers $M_i$ can be written as graphs over $S^n$, $M_i=\graph u_i$, and
the metric in $N$ can be expressed as

\begin{equation}
d\bar s^2={(dx^0)}^2+\s_{ij}(x^0,x) dx^i dx^j.
\end{equation}

\cvm
The 0-th component of the Gau{\ss} formula yields

\begin{equation}
v^{-1}h_{ij}=-u_{ij}+\bar h_{ij},
\end{equation}
\nd where
\begin{equation}
v^2=1+\s^{ij}u_iu_j,
\end{equation}

\cvm\nd
and  $(\bar h_{ij})$ is the second fundamental form of the level hypersurfaces
$\{x^0=\const\}$. 

\cvm
\br\lar{2.2}
It is well known that these level hypersurfaces are strictly
convex if $K_N\le 0$, and, hence, that there exists a strictly convex function
$\chi\in C^2(\bar\Om)$, cf. \ci[Remark 0.5]{cg4}.
\er

\cvb
\section{The evolution problem}\las 3

\cvb
Solving the problem \re{0.1} consists of two steps: first, one has to prove a
priori estimates, and secondly, one has to find a procedure which, with the help
of the priori estimates, leads to a solution of the problem.

\cvm
In Lorentzian manifolds the evolution method is the method of choice, but in
Riemannian manifolds this approach requires the sectional curvatures of the
ambient space to be non-positive. There is an alternative method---successive
approximation---but it is only applicable when the a priori estimates also apply
to the {\it elliptic regularizations} of the  curvature functions in mind, cf.
\ci{cg4}. Though the class $(K)$ is closed under elliptic regularization, see
\ci[Section 1]{cg4}, this is not valid for the subclass of functions satisfying the
additional property \re{0.3}. For that reason we require that a Riemannian space
has non-positive sectional curvature.

\cvb

We want to prove that the equation 

\begin{equation}
F=f
\end{equation}

\cvm\nd
has a solution. For technical reasons, it is convenient to solve instead the
equivalent equation

\begin{equation}\lae{3.2}
\F(F)=\F(f),
\end{equation}

\cvm\nd
where $\F$ is a real function defined on $\R[]_+$ such that

\begin{equation}
\dot\F>0\q \tup{and}\q \ddot\F\le 0.
\end{equation}

For notational reasons, let us abbreviate
\begin{equation}
\tilde f=\F(f).
\end{equation}

We also point out that we may---and shall---assume without loss of generality
that $F$ is homogeneous of degree 1 if $F$ is of class $(K)$.

\cvm

To solve \re{3.2} we look at the evolution problem

\begin{equation}\lae{3.5}
\begin{aligned}
\dot x&=-\s(\F-\tilde f)\n,\\
x(0)&=x_0,
\end{aligned}
\end{equation}

\cvm\nd
where $x_0$ is an embedding of an initial strictly convex, compact, space-like
hypersurface $M_0$, $\F=\F(F)$, and $F$ is evaluated at the principal curvatures
of the flow hypersurfaces $M(t)$, or, equivalently, we may assume that $F$
depends on the second fundamental form $(h_{ij})$ and the metric $(g_{ij})$ of
$M(t)$; $x(t)$ is the embedding of $M(t)$ and $\s$ the signature of the  normal 
$\n=\n(t)$---past directed, if N is Lorentzian, resp. the outward normal, if $N$ is
Riemannian.

\cvm

This is a parabolic problem, so short-time existence is guaranteed---the proof
in the Lorentzian case is identical to that in the Riemannian case, cf. \ci[p.
622]{cg2}---, and under suitable assumptions, we shall be able to prove that the
solution exists for all time and converges to a stationary solution if $t$ goes to
infinity.

\cvm
There is a slight ambiguity in the notation, since we also call the
evolution parameter \tit{time}, but this lapse shouldn't cause any
misunderstandings.

\cvm
Next, we want to show how the metric, the second fundamental form, and the
normal vector of the hypersurfaces $M(t)$ evolve. All time derivatives are
\tit{total} derivatives. The proofs are identical to those of the corresponding
results in a Riemannian setting, cf. \ci[Section 3]{cg2}, and will be omitted.

\cvm
\bl[Evolution of the metric]
The metric $g_{ij}$ of $M(t)$ satisfies the evolution equation
\begin{equation}
\dot g_{ij}=-2\s(\F-\tilde f)h_{ij}.
\end{equation}
\el

\cvm
\bl[Evolution of the normal]
The normal vector evolves according to
\begin{equation}\lae{3.7}
\dot \n=\nabla_M(\F-\tilde f)=g^{ij}(\F-\tilde f)_i x_j.
\end{equation}
\el

\cvm
\bl[Evolution of the second fundamental form]
The second fundamental form evolves according to
\begin{equation}\lae{3.8}
\dot h_i^j=(\F-\tilde f)_i^j+\s (\F-\tilde f) h_i^k h_k^j+\s (\F-\tilde f) \riema
\al\bet\ga\de\n^\al x_i^\bet \n^\ga x_k^\de g^{kj}
\end{equation}
and
\begin{equation}
\dot h_{ij}=(\F-\tilde f)_{ij}-\s (\F-\tilde f) h_i^k h_{kj}+\s (\F-\tilde f) \riema
\al\bet\ga\de\n^\al x_i^\bet \n^\ga x_j^\de.
\end{equation}
\el

\cvm
\bl[Evolution of $(\F-\tilde f)$]
The term $(\F-\tilde f)$ evolves according to the equation

\begin{multline}\lae{3.10}
{(\F-\tilde f)}^\prime-\dot\F F^{ij}(\F-\tilde f)_{ij}=
\s \dot \F
F^{ij}h_{ik}h_j^k (\F-\tilde f)\\+\s\tilde f_\al\n^\al (\F-\tilde f)
- \tilde f_{\n^\al}x^\al_i(\F- \tilde
f)_jg^{ij}\\
+\s\dot\F F^{ij}\riema \al\bet\ga\de\n^\al x_i^\bet \n^\ga x_j^\de
(\F-\tilde f),
\end{multline}

\cvm
\nd
where
\begin{equation}
(\F-\tilde f)^{\prime}=\frac{d}{dt}(\F-\tilde f)
\end{equation}
and
\begin{equation}
\dot\F=\frac{d}{dr}\F(r).
\end{equation}
\el

\cvm
From \re{3.8} we deduce with the help of the Ricci identities a parabolic equation
for the second fundamental form

\cvm
\bl
The mixed tensor $h_i^j$ satisfies the parabolic equation

\begin{equation}\lae{3.13}
\begin{aligned}
\dot h_i^j&-\dot\F F^{kl}h_{i;kl}^j\\[\cma]&=\s \dot\F
F^{kl}h_{rk}h_l^rh_i^j-\s\dot\F F h_{ri}h^{rj}+\s (\F-\tilde f)
h_i^kh_k^j\\[\cma] 
&\hp{+}-\tilde f_{\al\bet} x_i^\al x_k^\bet g^{kj}+\s \tilde f_\al\n^\al
h_i^j-\tilde f_{\al\n^\bet}(x^\al_i x^\bet_kh^{kj}+x^\al_l x^\bet_k h^{k}_i\, g^{lj})\\
&\hp{=}
-\tilde f_{\n^\al\n^\bet}x^\al_lx^\bet_kh^k_ih^{lj}-\tilde f_{\n^\bet} x^\bet_k
h^k_{i;l}\,g^{lj}  +\s\tilde f_{\n^\al}\n^\al h^k_i h^j_k\\
&\hp{=}+\dot\F
F^{kl,rs}h_{kl;i}h_{rs;}^{\hphantom{rs;}j}+2\dot \F F^{kl}\riema \al\bet\ga\de x_m^\al x_i ^\bet x_k^\ga
x_r^\de h_l^m g^{rj}\\
&\hp{=}-\dot\F F^{kl}\riema \al\bet\ga\de x_m^\al x_k ^\bet x_r^\ga x_l^\de
h_i^m g^{rj}-\dot\F F^{kl}\riema \al\bet\ga\de x_m^\al x_k ^\bet x_i^\ga x_l^\de h^{mj} \\
&\hp{=}+\s\dot\F F^{kl}\riema \al\bet\ga\de\n^\al x_k^\bet\n^\ga x_l^\de h_i^j-\s\dot\F F
\riema \al\bet\ga\de\n^\al x_i^\bet\n^\ga x_m^\de g^{mj}\\
&\hp{=}+\s (\F-\tilde f)\riema \al\bet\ga\de\n^\al x_i^\bet\n^\ga x_m^\de g^{mj}+\ddot \F
F_i F^j\\ &\hp{=}+\dot\F F^{kl}\bar R_{\al\bet\ga\de;\e}\{\n^\al x_k^\bet x_l^\ga x_i^\de
x_m^\e g^{mj}+\n^\al x_i^\bet x_k^\ga x_m^\de x_l^\e g^{mj}\}.
\end{aligned}
\end{equation}
\el

\cvm
The proof is identical to that of the corresponding result in \ci[Lemma 3.5]{cg8};
we only have to keep in mind that $f$ now also depends on the normal.

\cvm
If we had assumed $F$ to be homogeneous of degree $d_0$ instead of 1, then, we
would have to replace the explicit term $F$---occurring twice in the preceding
lemma---by $d_0F$.

\cvm
\br\lar{3.6}
In view of the maximum principle, we immediately deduce from \re{3.10} that the
term $(\F-\tilde f)$ has a sign during the evolution if it has one at the beginning,
i.e., if the starting hypersurface $M_0$ is the upper barrier $M_2$, then
$(\F-\tilde f)$ is non-negative, or equivalently,

\begin{equation}\lae{3.14}
F\ge f,
\end{equation}

\cvm\nd
while in case $M_0=M_1$, $(\F- \tilde f)$ is non-positive, or equivalently,

\begin{equation}\lae{3.15}
F\le   f.
\end{equation}
\er

\cvb
\section{Lower order estimates}\las{4}

\cvb
We consider the evolution problem \re{3.5} with $\F(r)=\log r$ and with initial
hypersurface $M_0=M_2$ if $N$ is Lorentzian resp. $M_0=M_1$ if $N$ is
Riemannian. Solutions exist in a maximal time interval $[0,T^*)$, $0<T^*\le \un$,
as long as the flow hypersurfaces stay in $\bar\Om$ and are smooth and strictly
convex.

\cvm
Let us first consider the Lorentzian case in more detail.

\subsection{The Lorentzian case} As we have already mentioned, the barriers
$M_i$ are then graphs over the compact Cauchy hypersurface $\so$ and this is
also valid for the flow hypersurfaces $M(t)$, $M(t)=\graph u(t)$.

\cvm
The scalar version of \re{3.5} is

\begin{equation}
\pde ut=-e^{-\psi} v (\F - \tilde f),
\end{equation}

\cvm
\nd where

\begin{equation}
v= \tilde v^{-1}=1-\abs{Du}^2.
\end{equation}

\cvm
As we have shown in \ci[Section 4]{cg8}, the flow hypersurfaces stay in $\bar\Om$
and are uniformly space-like, i.e. the term $ \tilde v$ is uniformly bounded.
Moreover,  $\tilde v$ satisfies a
useful parabolic equation that we shall exploit to estimate the principal
curvatures of the hypersurfaces $M(t)$ from above.

\cvm
\bl[Evolution of $\tilde v$]
Consider the flow \re{3.5} in the distinguished coordinate system associated
with $\so$. Then, $\tilde v$ satisfies the evolution equation

\begin{equation}\lae{4.3}
\begin{aligned}
\dot{\tilde v}-\dot\F F^{ij}\tilde v_{ij}=&-\dot\F F^{ij}h_{ik}h_j^k\tilde v
+[(\F-\tilde f)-\dot\F F]\h_{\al\bet}\n^\al\n^\bet\\
&-2\dot\F F^{ij}h_j^k x_i^\al x_k^\bet \h_{\al\bet}-\dot\F F^{ij}\h_{\al\bet\ga}x_i^\bet
x_j^\ga\n^\al\\
&-\dot\F F^{ij}\riema \al\bet\ga\de\n^\al x_i^\bet x_k^\ga x_j^\de\h_\e x_l^\e g^{kl}\\
&-\tilde f_\bet x_i^\bet x_k^\al \h_\al g^{ik}- \tilde f_{\n^\bet}x^\bet_k h^{ik}x^\al_i\h_\al,
\end{aligned}
\end{equation}

\cvm\nd
where $\h$ is the covariant vector field $(\h_\al)=e^{\psi}(-1,0,\dotsc,0)$.
\el

\cvm
For a proof see \ci[Lemma 4.4]{cg8}; we only have to keep in mind that, now, $f$ 
also depends on the normal.

\cvm
\bc
Let $ \tilde \f=e^{\lam \tilde v}$, then, $ \tilde \f$ satisfies the evolution inequality

\begin{equation}\lae{4.4}
\begin{aligned}
\dot {\tilde \f}-\dot\F F^{ij} \tilde \f_{ij}\le &-\tfrac\lam2\dot\F F^{ij}h_{ik}h_j^k
e^{\lam \tilde v}- \tfrac{\lam^2}2 \dot\F F^{ij} \tilde v_i \tilde v_j e^{\lam \tilde
v}\\[\cma]
 &+c\lam \dot\F F^{ij}g_{ij} e^{\lam \tilde v} + c[(\F- \tilde f)+1]e^{\lam \tilde
v}\\[\cma]
 &+c\lam \dot\F F^{ij}g_{ij} e^{\lam \tilde v} + c(\dot\F)^{-1} \tilde F^{ij}g_{ij} e^{\lam
\tilde v},
\end{aligned}
\end{equation}

\cvm
\nd
where $( \tilde F^{ij})=(F^{ij})^{-1}$,  $c$ is a known constant, and where we also
used the estimate \re{3.14}.
\ec

\bp
We have

\begin{equation}
\begin{aligned}
\dot {\tilde \f}-\dot\F F^{ij} \tilde \f_{ij}=[\dot{ \tilde v}-\dot\F F^{ij} \tilde
v_{ij}]\lam e^{\lam \tilde v}- \lam^2 \dot\F F^{ij} \tilde v_i \tilde v_j e^{\lam \tilde v}.
\end{aligned}
\end{equation}

\cvm
The non-trivial terms in \re{4.3} are estimated as follows

\begin{equation}
-2\dot\F F^{ij}h_j^k x_i^\al x_k^\bet \h_{\al\bet}\lam e^{\lam \tilde v}\le \tfrac\lam2\dot\F
F^{ij}h_{ik}h_j^k e^{\lam \tilde v}+c\lam \dot\F F^{ij}g_{ij} e^{\lam \tilde v},
\end{equation}

\cvm
\nd and
\begin{equation}\lae{4.7}
\begin{aligned}
\abs{\tilde f_{\n^\bet}x^\bet_k h^{ik}x^\al_i\h_\al}\lam e^{\lam \tilde v}&\le  \,\abs{ \tilde
f_{\n^\bet} \tilde v_kx^\bet_lg^{kl}}\lam e^{\lam \tilde v}+ c\nnorm{ \tilde f_{\n^\bet}}\lam
e^{\lam \tilde v}\\[\cma]
&\le c\nnorm{ \tilde f_{\n^\bet}}\lam e^{\lam \tilde v}+\tfrac{\lam^2}2 \dot\F F^{ij}  \tilde
v_i
\tilde v_j e^{\lam \tilde v}\\
&\qq\qq\qq\,\,+ c(\dot\F)^{-1} \tilde F^{ij}g_{ij} e^{\lam
\tilde v}.
\end{aligned}
\end{equation}

\cvm
\nd With the help these estimates inequality \re{4.4} is easily derived; in the
first inequality of \re{4.7} we used

\begin{equation}\lae{4.8}
\tilde v_i=\h_{\al\bet}x_i^\bet\n^\al+\h_\al\n_i^\al
\end{equation}

\cvm
\nd together with the Weingarten equation.
\ep

\subsection{The Riemannian case} As we have shown in \ci[Sections 5 \& 6]{cg2},
the flow hypersurfaces can be written as graphs over a geodesic unit sphere,
$M(t)=\graph u(t)$. The scalar version of \re{3.5} now looks like

\begin{equation}
\pde ut=-(\F- \tilde f)v.
\end{equation}

\cvm
Moreover, all flow hypersurfaces stay in $\bar\Om$ and $v$ is uniformly bounded
in view of the convexity of the $M(t)$.

\cvm
\bl[Evolution of $v$]
The quantity $v$ satisfies the parabolic equation

\begin{equation}\lae{4.10}
\begin{aligned}
\dot{ v}-\dot\F F^{ij} v_{ij}=&-\dot\F F^{ij}h_{ik}h_j^k v-2 v^{-1} \dot\F
F^{ij}v_iv_j \\
&+[(\F- f)-\dot\F F]\h_{\al\bet}\n^\al\n^\bet v^2\\
&+2\dot\F F^{ij}h_j^k x_i^\al x_k^\bet \h_{\al\bet} v^2+\dot\F F^{ij}\h_{\al\bet\ga}x_i^\bet
x_j^\ga\n^\al v^2\\
&+\dot\F F^{ij}\riema \al\bet\ga\de\n^\al x_i^\bet x_k^\ga x_j^\de\h_\e x_l^\e g^{kl} v^2\\
&+\tilde f_\bet x_i^\bet x_k^\al \h_\al g^{ik}v^2+\tilde f_{\n^\bet}x^\bet_k h^{ik}x^\al_i
\h_\al v^2,
\end{aligned}
\end{equation}

\cvm\nd
where $\h$ is the covariant vector field $(\h_\al)=(1,0,\dotsc,0)$.
\el

\cvm
For a proof see \ci[Lemma 7.3]{cg2}.

\cvm
Similar as in the Lorentzian case we obtain a parabolic inequality for $e^{\lam v}$.

\cvm
\bc
Let $  \f=e^{\lam  v}$, then, $  \f$ satisfies the evolution inequality

\begin{equation}\lae{4.11}
\begin{aligned}
\dot { \f}-\dot\F F^{ij}  \f_{ij}\le &-\tfrac\lam2\dot\F F^{ij}h_{ik}h_j^k
e^{\lam  v}- \tfrac{\lam^2}2 \dot\F F^{ij}  v_i  v_j e^{\lam 
v}\\[\cma]
 &+c\lam \dot\F F^{ij}g_{ij} e^{\lam  v} + c[-(\F-  f)+1]e^{\lam 
v}\\[\cma]
 &+c\lam \dot\F F^{ij}g_{ij} e^{\lam  v} + c(\dot\F)^{-1} \tilde  F^{ij}g_{ij} e^{\lam
 v},
\end{aligned}
\end{equation}

\cvm
\nd
where $( \tilde F^{ij})=(F^{ij})^{-1}$, $c$ is a known constant, and where we also
used the estimate \re{3.15}.
\ec

\cvb
\section{$C^2$- estimates in Lorentzian space}\las{5}

\cvb
Let $M(t)$ be a solution of the evolution problem \re{3.5} with initial
hypersurface $M_0=M_2$, defined on a maximal time interval $I=[0,T^*)$. We
assume that $F$ is of class $(K)$, homogeneous of
degree 1, and satisfies the condition \re{0.3}; we choose
$\F(r)=\log r$. 

\cvm
Furthermore, we suppose that there exists
a strictly convex function $\chi\in C^2(\bar \Om)$, i.e. there holds

\begin{equation}
\chi_{\al\bet}\ge c_0\bar g_{\al\bet}
\end{equation}

\cvm\nd
with a positive constant $c_0$.

\cvm
We observe that

\begin{equation}\lae{5.2}
\begin{aligned}
\dot\chi-\dot\F F^{ij}\chi_{ij}&=[(\F-\tilde f)-\dot\F F]\chi_\al\n^\al-\dot\F
F^{ij}\chi_{\al\bet} x_i^\al x_j^\bet\\
&\le [(\F-\tilde f)-\dot\F F] \chi_\al\n^\al-c_0\dot\F F^{ij} g_{ij},
\end{aligned}
\end{equation}

\cvm\nd
where we used the homogeneity of $F$.

\cvm
From \rr{3.6} we infer 

\begin{equation}\lae{5.4}
\F\ge \tilde f\qq \tup{or}\qq F\ge f,
\end{equation}

\cvm\nd
and from the results in \rs{4} that the flow stays in the compact set $\bar \Om$,
and that $ \tilde v$ is uniformly bounded.

\cvm
We are now able to prove

\cvm
\bl\lal{5.1}
Let $F$ be of class $(K)$ satisfying \re{0.3}. Then, the principal curvatures of the
evolution hypersurfaces $M(t)$ are uniformly bounded.
\el

\bp
Let $\f$ and $w$ be defined respectively by

\begin{align}
\f&=\sup\set{{h_{ij}\h^i\h^j}}{{\norm\h=1}},\lae{5.4b}\\
w&=\log\f+e^{\lam \tilde v}+\m\chi,\lae{5.6}
\end{align}

\cvm\nd
where $\lam,\m$ are large positive parameters to be specified later. We claim that
$w$ is bounded for a suitable choice of $\lam,\m$.

\cvm
Let $0<T<T^*$, and $x_0=x_0(t_0)$, with $ 0<t_0\le T$, be a point in $M(t_0)$ such
that

\begin{equation}
\sup_{M_0}w<\sup\set {\sup_{M(t)} w}{0<t\le T}=w(x_0).
\end{equation}

We then introduce a Riemannian normal coordinate system $(\x^i)$ at $x_0\in
M(t_0)$ such that at $x_0=x(t_0,\x_0)$ we have

\begin{equation}
g_{ij}=\de_{ij}\q \tup{and}\q \f=h_n^n.
\end{equation}

\cvm
Let $\tilde \h=(\tilde \h^i)$ be the contravariant vector field defined by

\begin{equation}
\tilde \h=(0,\dotsc,0,1),
\end{equation}
and set
\begin{equation}
\tilde \f=\frac{h_{ij}\tilde \h^i\tilde \h^j}{g_{ij}\tilde \h^i\tilde \h^j}\raise 2pt
\hbox{.}
\end{equation}

\cvm
$\tilde \f$ is well defined in neighbourhood of $(t_0,\x_0)$.

\cvm
Now, define $\tilde w$ by replacing $\f$ by $\tilde \f$ in \re{5.6}; then, $\tilde w$
assumes its maximum at $(t_0,\x_0)$. Moreover, at $(t_0,\x_0)$ we have

\begin{equation}
\dot{\tilde \f}=\dot h_n^n,
\end{equation}

\cvm\nd
and the spatial derivatives do also coincide; in short, at $(t_0,\x_0)$ $\tilde \f$
satisfies the same differential equation \re{3.13} as $h_n^n$. For the sake of
greater clarity, let us therefore treat $h_n^n$ like a scalar and pretend that $w$
is defined by

\begin{equation}
w=\log h_n^n+ e^{\lam\tilde v}+\m\chi.
\end{equation}

\cvm
At $(t_0,\x_0)$ we have $\dot w\ge 0$, and, in view of the maximum principle, we
deduce from \re{0.3}, \re{1.17}, \re{3.13}, \re{4.4}, and \re{5.2}

\begin{equation}\lae{5.13}
\begin{aligned}
0\le&\msp[3] \dot\F F h_n^n-(\F-\tilde f) h_n^n+\lam c_1e^{\lam \tilde v}
-\tfrac\lam 2\e_0\dot\F F H
e^{\lam \tilde v}\\ 
& +\lam c_1\dot\F F^{ij} g_{ij}e^{\lam \tilde v}+c_1  H e^{\lam
\tilde v}\\ &+ (\lam e^{\lam \tilde v}+\m) c_1 [(\F-\tilde f)+\dot\F F]-\m c_0 \dot\F
F^{ij} g_{ij}\\ &+\dot\F F^{ij}(\log h_n^n)_i(\log h_n^n)_j\\
&+\{\ddot\F F_n F^n +\dot\F
F^{kl,rs}h_{kl;n}h_{rs;}^{\hphantom{rs;}n}\}(h_n^n)^{-1},
\end{aligned}
\end{equation}

\cvm\nd
where we have estimated bounded terms by a constant $c_1$, assumed that
$h_n^n, \lam$, and $\m$ are larger than $1$, and used \re{5.4}. 

\cvm
Now, the last term in \re{5.13} is estimated from above by

\begin{equation}\lae{5.15}
\{\ddot\F F_n F^n+\dot\F F^{-1} F_n F^n\}(h_n^n)^{-1}-\dot \F F^{ij}
h_{in;n}h_{jn;}^{\hphantom{jn;}n}(h_n^n)^{-2},
\end{equation}

\cvm\nd
cf. \re{1.10}, where  the sum in the braces vanishes,
due to the choice of
$\F$. Moreover, because of the Codazzi equation, we have

\begin{equation}
h_{in;n}=h_{nn;i}+\riema \al\bet\ga\de\n^\al x_n^\bet x_i^\ga x_n^\de,
\end{equation}

\cvm\nd
and hence, using the  abbreviation $\bar R_i$ for the curvature term, we conclude
that \re{5.15} is bounded from above by

\begin{equation}
-(h_n^n)^{-2}\dot\F F^{ij}(h_{n;i}^n+\bar R_i)(h_{n;j}^n+\bar R_j).
\end{equation}

\cvm
Thus, the terms in \re{5.13} containing the derivatives of $h_n^n$ are estimated
from above by

\begin{equation}
-2\dot\F F^{ij}(\log h_n^n)_i\bar R_j(h_n^n)^{-1}.
\end{equation}

\cvm
Moreover,  $Dw$ vanishes at $\x_0$, i.e.

\begin{equation}
D\log h_n^n=-\lam e^{\lam \tilde v} D\tilde v-\m D\chi,
\end{equation}

\cvm\nd
where only $D\tilde v$ deserves further consideration.

\cvm
Replacing then $D\tilde v$ by the right-hand side of \re{4.8}, and using the
Weingarten equation as well as the simple
observation

\begin{equation}
\abs{F^{ij}h_j^k\h_k}\le \norm\h F
\end{equation}

\cvm\nd
for any vector field $(\h_k)$, cf. \ci[Lemma 7.4]{cg2}, we finally conclude from
\re{5.13}

\begin{equation}
\begin{aligned}
0\le &\msp[3]\dot\F F h_n^n-(\F-\tilde f)h_n^n+\lam c_1e^{\lam \tilde v} + c_1 H e^{\lam
\tilde v}-
\tfrac\lam 2\e_0\dot\F F H
e^{\lam \tilde v}\\
&+(\lam e^{\lam \tilde v}+\m)c_1[(\F-\tilde f)+\dot\F F]+\lam c_1 e^{\lam \tilde v}\dot\F
F^{ij} g_{ij}\\ &-\m [c_0-c_1 (h_n^n)^{-1}]\dot\F F^{ij} g_{ij}
\end{aligned}
\end{equation}

\cvm
Then, if we suppose $h_n^n$ to be so large that

\begin{equation}
c_1\le \tfrac{1}{2}c_0 h_n^n,
\end{equation}

\cvm\nd
and if we choose $\lam, \m$ such that

\begin{align}
4&\le \lam\e_0\\
\intertext{and}
8\lam c_1&\le \m c_0
\end{align}

\cvm\nd
we derive

\begin{equation}\lae{5.22}
\begin{aligned}
0\le&\msp[2] -\tfrac{1}{4}\lam\e_0 \dot\F F H e^{\lam \tilde v}-(\F-\tilde f) h_n^n+
c_1 H e^{\lam \tilde v}\\ &+(\lam e^{\lam \tilde v}+\m)c_1 [(\F-\tilde f)+\dot\F F].
\end{aligned}
\end{equation}

\cvm
We now observe that $\dot \F F=1$, and deduce in view of \re{5.4} that $h_n^n$ is
a priori bounded at $(t_0,\x_0)$.
\ep

The result of \rl{5.1} can be restated as a uniform estimate for the functions
$u(t)\in C^2(\so)$. Since, moreover, the principal curvatures of the flow
hypersurfaces are not only bounded, but also uniformly bounded away from zero,
in view of \re{5.4} and the assumption that $F$ vanishes on $\pa \C_+$, we
conclude that $F$ is uniformly elliptic on $M(t)$.

\cvb
\section{$C^2$- estimates in Riemannian space}\las 6

\cvb
If $N$ is Riemannian with $K_N\le0$, we use $M_1$ as the initial hypersurface for
the evolution. This is the only change in the settings with regard to the
Lorentzian case. Due to this choice we now have

\begin{equation}\lae{6.1}
\F- \tilde f\le 0
\end{equation}

\cvm\nd
during the evolution, cf. \rr{3.6}..

\cvm
For the $C^2$- estimates we have to prove upper bounds for the principal
curvatures as well as a strictly positive lower bound for $F$

\begin{equation}\lae{6.2}
0<\e_2\le F,
\end{equation}

\cvm\nd
since, in view of \rr{3.6}, we presently only know that \re{3.15} is valid.

\cvm
Furthermore, as we have seen in \rr{2.2}, there exists a strictly convex function
$\chi$.

\cvm
We now prove the corresponding result to \rl{5.1}

\cvm
\bl\lal{6.1}
Let $F\in (K)$ such that \re{0.3} is satisfied, and let $M(t)$ be the solutions of
\re{3.5} with initial hypersurface $M_0=M_1$. Then, the principal curvatures of
the $M(t)$ are uniformly bounded from above.
\el

\bp
The proof is identical to that of \rl{5.1}. We define $\f$ and $w$ as in \re{5.4b} and
\re{5.6}, where of course $ \tilde v$ is replaced by $v$, and apply the maximum
principle to $w$.

\cvm
In a point where the maximum principle is applied we obtain in view of \re{1.17}
an inequality that corresponds to the similar inequality \re{5.22}

\begin{equation}\lae{6.3}
\begin{aligned}
0\le &-\tfrac\lam 4\e_0\dot \F F H e^{\lam v} + (\F- \tilde f) h^n_n +c_1 H e^{\lam v}\\
& +(\lam e^{\lam v}+\m) c_1 [-(\F- \tilde f) +\dot\F F].
\end{aligned}
\end{equation}

\cvm
In deriving this inequality, we also used the simple estimate

\begin{equation}
\dot\F F^{ij}h_{ik}h_j^k\le \dot \F F h^n_n.
\end{equation}

\cvm
From \re{6.3} we immediately get the required a priori estimate because of
\re{6.1}.
\ep

\cvm
\bl
Assume that $K_N\le 0$, then, there is a positive constant $\e_2$ such that the
estimate \re{6.2} is valid.
\el

\bp
We proceed similar as in the proof of \ci[Lemma 8.3]{cg2}. Consider the function

\begin{equation}
w=-(\F- \tilde f) +\m\chi,
\end{equation}

\cvm
\nd where $\m$ is large. Let $0<T< T^*$ and suppose

\begin{equation}
\sup_{M_0}w<\sup\set{\sup_{M(t)} w}{0\le t\le T}.
\end{equation}

\cvm
Then, there exists $x_0=x_0(t_0)\in \so$, $0<t_0\le T$, such that

\begin{equation}
w(x_0)=\sup\set{\sup_{M(t)} w}{0\le t\le T}.
\end{equation}

\cvm
From \re{3.10}, \re{5.2}, and the maximum principle we then infer

\begin{multline}\lae{6.8}
0\le
-\dot \F
F^{ij}h_{ik}h_j^k (\F-\tilde f)-\tilde f_\al\n^\al (\F-\tilde f)
+ \tilde f_{\n^\al}x^\al_i(\F- \tilde
f)_jg^{ij}\\
-\dot\F F^{ij}\riema \al\bet\ga\de\n^\al x_i^\bet \n^\ga x_j^\de
(\F-\tilde f) +\m c_1 [1-(\F - \tilde f)]-\m c_0 \dot\F F^{ij}g_{ij}.
\end{multline}

\cvm
Let $\ka$ be an upper bound for the principal curvatures, then, the first term on
the right-hand side can be estimated from above by

\begin{equation}
-\dot \F F^{ij}h_{ik}h_j^k (\F-\tilde f)\le -\dot\F F\ka (\F- \tilde f).
\end{equation}

\cvm
The term involving the Riemann curvature tensor is non-positive since $K_N\le
0$, and, hence, we deduce

\begin{equation}\lae{6.10}
0\le \m c_1[1-(\F- \tilde f)] -\m c_0 \dot\F F^{ij}g_{ij}
\end{equation}

\cvm\nd
for $\m\ge 1$, and we obtain an a priori estimate for $-(\F- \tilde f)$, since

\begin{equation}
F^{ij}g_{ij}\ge F(1,\dotsc, 1)
\end{equation}

\cvm
\nd and $\dot \F=F^{-1}$ is the dominating term in \re{6.10}.
\ep

\cvb
\br
The assumption $K_N\le 0$ was only necessary to obtain a uniform bound for the
principal curvatures during the {\it evolution}. For stationary solutions

\begin{equation}\lae{6.12}
\fmo M= f(x,\n)
\end{equation}

\cvm
\nd the proof of \rl{6.1} would yield a priori estimates for the principal
curvatures in arbitrary Riemannian manifolds as long as $M$ is a graph in a
Gaussian coordinate system, lower order estimates are valid, and there exists a
strictly convex function in a neighbourhood of $M$.

\cvm
This could be used to solve the Dirichlet problem for the equation \re{6.12},
since in the existence proof for Dirichlet problems a deformation process is
used instead of an evolutionary approximation, cf. \ci{os}.
\er

\cvb
\section{Convergence to a stationary solution}\las 7

We only consider the Lorentzian case since the essential arguments do not
depend on the nature of the ambient space.  Let
$M(t)$ be the flow with initial hypersurface $M_0=M_2$. Let us look at the scalar
version of the flow
\re{3.5}

\begin{equation}\lae{7.1}
\pde ut=-e^{-\psi}v(\F-\tilde f).
\end{equation}

\cvm\nd
This is  a scalar parabolic differential equation defined on the cylinder

\begin{equation}
Q_{T^*}=[0,T^*)\times \so
\end{equation}

\cvm\nd
with initial value $u(0)=u_2\in C^{4,\al}(\so)$. In view of the a priori estimates,
which we have established in the preceding sections, we know that

\begin{equation}
{\abs u}_\low{2,0,\so}\le c
\end{equation}
and
\begin{equation}
\F(F)\,\tup{is uniformly elliptic in}\,u
\end{equation}

\cvm\nd
independent of $t$. Moreover, $\F(F)$ is concave, and thus, we can apply
the regularity results of \ci[Chapter 5.5]{nk} to conclude that uniform
$C^{2,\al}$-estimates are valid, leading further to uniform $C^{4,\al}$-estimates
due to the regularity results for linear operators.

Therefore, the maximal time interval is unbounded, i.e. $T^*=\un$.

\cvm
Now, integrating \re{6.1} with respect to $t$, and observing that the right-hand
side is non-positive, yields

\begin{equation}
u(0,x)-u(t,x)=\int_0^te^{-\psi}v(\F-\tilde f)\ge c\int_0^t(\F-\tilde f),
\end{equation}
i.e.,
\begin{equation}
\int_0^\un \abs{\F-\tilde f}<\un\qq\A\msp x\in \so
\end{equation}

\cvm\nd
Hence, for any $x\in\so$ there is a sequence $t_k\rightarrow \un$ such that
$(\F-\tilde f)\rightarrow 0$.

\cvm
On the other hand, $u(\cdot,x)$ is monotone decreasing and therefore

\begin{equation}
\lim_{t\rightarrow \un}u(t,x)=\tilde u(x)
\end{equation}

\cvm\nd
exists and is of class $C^{4,\al}(\so)$ in view of the a priori estimates. We, finally,
conclude that $\tilde u$ is a stationary solution of our problem, and that

\begin{equation}
\lim_{t\rightarrow \un}(\F-\tilde f)=0.
\end{equation}

\cvb

\end{document}